\documentclass{article}
\usepackage{amssymb}
\usepackage{amsmath}

\input{tcilatex}

\begin{document}

\section{Pre-Hausdorff Spaces}

\bigskip

JAY STINE

\textit{Department of Mathematics, Misericordia University, Dallas, PA,
18612, U.S.A.}

\bigskip

M. V. MIELKE

\textit{Department of Mathematics, University of Miami, Coral Gables, FL,
33124, U.S.A.}

\bigskip

\textbf{Abstract.} \ This paper introduces three separation conditions for
topological spaces, called T$_{0,1}$, T$_{0,2}$ ("pre-Hausdorff"), and T$%
_{1,2}$. \ These conditions generalize the classical T$_{1}$ and T$_{2}$
separation axioms, and they have advantages over them topologically which we
discuss. \ We establish several different characterizations of pre-Hausdorff
spaces, and a characterization of Hausdorff spaces in terms of
pre-Hausdorff. \ We also discuss some classical Theorems of general topology
which can or cannot be generalized by replacing the Hausdorff condition by
pre-Hausdorff.

\bigskip

\textbf{Mathematics Subject Classifications (2000):} \ 18B30, 54A05, 54D10.

\bigskip

\textbf{Keywords:} \ topological separation properties, topological
category, left adjoint, sober space

\bigskip

\textbf{Introduction}

\bigskip

The notion of separation is fundamental to topology. \ Even so, the
classical separation axioms (T$_{0}$, ..., T$_{4}$) are sometimes overlooked
in a first course or, alternatively, some consider the T$_{2}$ axiom (for
instance) as being sufficiently weak that all spaces are assumed T$_{2}$ and
no further consideration is given to separation. \ While this may be
reasonable in some settings, it is certainly not in others. \ Analysis often
takes place in the setting of metric spaces, which are T$_{2},$whereas
geometry often uses pseudometric (more generally, uniform) spaces which are
not necessarily T$_{2}$. \ Herrlich argues (in $\left[ 7\right] $) that
"there are sufficient reasons for topologists to pay serious attention to
non-Hausdorff spaces ... finite Hausdorff spaces are rare and not very
interesting ... a 14-element set carries just a single Hausdorff topology
but 98,484,324,257,128,207,032,183 T$_{0}$ topologies". \ In this paper we
define a generalized Hausdorff separation condition called pre-Hausdorff,
which is satisfied by many important non-Hausdorff spaces. \ In $\left[ 15%
\right] $ it is shown that a uniform space is T$_{0}$ if and only if it is T$%
_{2}$, and the proof of this reveals that all uniform spaces are
pre-Hausdorff. \ Following Herrlich's argument above, the worthiness of
studying pre-Hausdorff spaces can be justified by their abundance: a
14-element set carries 190,899,322 distinct pre-Hausdorff topologies (see
Corollary 2.3).

\ The paper is organized as follows. \ Section 1 contains definitions of
three separation axioms for topological spaces and examples to show how they
are related. \ We prove that the categories formed by the spaces which
satisfy these axioms are topological categories and, further, that these
categories are reflective in the category of topological spaces. \ In
section 2, we give several characterizations of pre-Hausdorff spaces in
terms of Hausdorff separation and some equivalence relations. \ Finally, in
section 3, we consider some classical Theorems of general topology which can
or cannot be generalized by replacing the Hausdorff condition by
pre-Hausdorff.

\ Throughout the paper, TOP will be used to denote the category of
topological spaces and continuous functions. \ For $i$ = 0,1,2, T$_{i}$-TOP
will denote the full subcategories of TOP consisting of the T$_{i}$ spaces
(see $\left[ 4\right] $, page 138)..

\bigskip

\textbf{1. T}$_{i,j}$ \textbf{- Spaces}

\bigskip

DEFINITION 1.1. A topological space X is called a T$_{i,j}$ - space (for 0 $%
\leq $ $i$ $<$ $j$ $\leq $ $2$) if and only if each pair of points a, b $\in 
$ X which has a T$_{i}$ - separation in X also has a T$_{j}$ - separation in
X.

\bigskip

NOTATION: The categories consisting of the T$_{i,j}$ - spaces, along with
continuous functions, will be denoted T$_{i,j}$ - TOP.

\bigskip

EXAMPLE 1.2. T$_{0,2}$ spaces have been refered to as pre-Hausdorff spaces
in the literature (see$\left[ 14\right] $).

\bigskip

EXAMPLE 1.3. T$_{0,1}$ spaces have been refered to as R$_{0}$ - spaces in
the literature (see$\left[ 8\right] $). An R$_{0}$ - space is a topological
space X which satisfies: x $\in $ $\overline{\left\{ \text{y}\right\} }$
(the topological closure of $\left\{ \text{y}\right\} $) if and only if y $%
\in $ $\overline{\left\{ \text{x}\right\} }$, for all pairs of points x, y $%
\in $ X. Evidently then, every neighborhood of x contains y if and only if
every neighborhood of y contains x. Now if x and y have no T$_{0}$
separation, then this condition is satisfied. On the other hand, if x and y
do have a T$_{0}$ separation, say x has a neighborhood not containing y,
then y must have a neighborhood not containing x; i.e., x and y must have a T%
$_{1}$ separation. Thus R$_{0}$ - spaces are exactly the T$_{0,1}$ - spaces.

\bigskip

EXAMPLE 1.4. Clearly T$_{j}$ spaces are always T$_{i,j}$, but T$_{i}$ spaces
need not be T$_{i,j}$. A Sierpinski space (i.e., a two-point set , say X = $%
\left\{ \text{0,1}\right\} $, with one proper open set, say $\left\{ \text{1}%
\right\} $), for instance, is T$_{\text{0}}$, but neither T$_{0,1}$ nor T$%
_{0,2}$; while a T$_{1}$ space which is not T$_{2}$ will not be T$_{1,2}$.
Furthermore, a T$_{i,,j}$ space need not be either T$_{i}$ or T$_{j}$. An
indiscrete space with more than one element, for instance, is T$_{i,,j}$ for
each $i,j$, but is not T$_{0}$.

\bigskip

EXAMPLE 1.5. Clearly T$_{0,2}$ spaces are both T$_{0,1}$ and T$_{1,2}$.
However, a T$_{1,2}$ space need not be T$_{0,1}$ (and, consequently, not T$%
_{0,2}$ either), as in the case of a Sierpinski space, for example.
Furthermore, a T$_{0,1}$ space may be neither T$_{0,2}$ nor T$_{1,2}.$ This
is the case if, for example, a space is T$_{1}$ but not T$_{2}$.

\bigskip

The following Theorem shows that the categories T$_{i,j}$ - TOP have a
desirable property that is not shared by the categories T$_{i}$ - TOP.

\bigskip

THEOREM 1.6. \textit{The full subcategories T}$_{i,j}$ \textit{- TOP are
themselves topological over SET (the category of sets and functions).
Moreover, their inclusions into TOP preserve initial lifts and, consequently
they preserve all limits.}

\textit{Proof:} We prove the Theorem for T$_{0,2}$ -TOP, the cases T$_{0,1}$
- TOP and T$_{1,2}$ - TOP being similar. Clearly the restriction of the
forgetful functor U : TOP $\rightarrow $ SET is both concrete and has
set-theoretic fibers. So we show that the structure induced on a set from an
arbitrary family of T$_{0,2}$ spaces yields a T$_{0,2}$ space. This will
also show that initial lifts in T$_{0,2}$ - TOP are computed as they are in
TOP and, thus, the inclusion functor preserves them. Suppose that (X, $\tau $%
) is the induced topological space on a set X from a family $\left\{ (\text{X%
}_{i}\text{, }\tau _{i}\text{)}\right\} _{i\in I}$ of T$_{0,2}$ spaces via a
family of functions $\left\{ \text{f}_{i}\text{ : X }\rightarrow \text{ X}%
_{i}\right\} _{i\in I}$. Further suppose that x, y $\in $ X have a T$_{0}$ -
separation in $\tau $ by, say, U$_{x}$ $\in $ $\tau $, where x $\in $ U$_{x}$
and y $\notin $ U$_{x}$. We can assume that U$_{x}$ is a basis element of $%
\tau $ so that U$_{x}$ = $\overset{n}{\underset{j=1}{\bigcap }}$ f$_{i_{j}}^{%
\text{-1}}$(V$_{j}$), where each V$_{j}$ is open in X$_{i_{j}}$ for each $j$
= 1, 2, ..., $n$. Then $\exists k$, 1 $\leq k\leq n$, with f$_{i_{k}}$(x) $%
\in $ V$_{k}$ and f$_{i_{k}}$(y) $\notin $ V$_{k}$; i.e., f$_{i_{k}}$(x) and
f $_{i_{k}}$(y) have a T$_{0}$ - separation in X$_{i_{k}}$. Since X$_{i_{k}}$
is T$_{0,2}$, $\exists $ nbhds. U and W of f$_{i_{k}}$(x) and f$_{i_{k}}$(y)
(resp.) such that U $\cap $ W = $\emptyset $. Therefore (X, $\tau $) is T$%
_{0,2}$.

\bigskip

COROLLARY 1.7. \ The inclusion functors \textit{inc}$_{i,j}$: T$_{i,j}$ -
TOP $\rightarrow $ TOP each have a left adjoint L$_{i,j}$.

\textit{Proof: }Note that any indiscrete space with two elements forms a
small (one element) cogenerating set for any of the categories T$_{i,j}$ -
TOP. \ Since the functors \textit{inc}$_{i,j}$ are continuous by Theorem
1.6, the result follows immediately from the Corollary on page 126 of $\left[
12\right] $.

\bigskip

Note: An explicit description of the left adjoint to \textit{inc}$_{0,2}$ is
given below, in the discussion following Theorem 2.18. \ Another description
of L$_{0,2}$ using transfinite recursion can be found in $\left[ 19\right] $%
, where this approach is then adapted to give an explicit description of the
functor L$_{0,1}:$ TOP $\rightarrow $ T$_{0,1}$ - TOP. \ It is also shown
there that these left adjoints are retractions.

\bigskip

\textbf{2. Pre-Hausdorff Spaces}

\bigskip

This section is concerned specifically with T$_{0,2}$-TOP, the category of
pre-Hausdorff spaces. \ In $\left[ 18\right] $ Steiner defines principal
topologies in terms of ultratopologies, and proves that a topological space
is principal if and only if arbitrary intersections of open sets are open
(such spaces are also refered to as Alexandroff spaces in the literature,
see $\left[ 2\right] ;$ however the term Alexandroff space also appears in a
different context, see $\left[ 3\right] $ ). \ The following result will be
used to gain insight into finite pre-Hausdorff spaces, and to count the
number of distinct pre-Hausdorff topologies on a a finite set. \ We note
that in this result, as in $\left[ 11\right] $ page 113, we make a
distinction between regular spaces and T$_{3}$ spaces; namely, that regular
spaces need not have closed points.

\bigskip

THEOREM 2.1. \textit{Suppose (X, }$\tau $) \textit{is a principal space (X, }%
$\tau $). \ \textit{The following are equivalent:}

\ \ (i) \textit{(X}, $\tau $\textit{)} \textit{is pre-Hausdorff.}

\ \ (ii) \ \textit{(X}, $\tau $\textit{)} \textit{is regular.}

\ \ (iii) \ \textit{(X, }$\tau $\textit{) has dimension 0; i.e., has a basis
consisting of clopen sets (see }$\left[ 9\right] $\textit{, page 10, B).}

\ (iv) \textit{The topos of sheaves on X is Boolean; i.e., the negation
operator }$\lnot :\tau \rightarrow \tau $ \textit{satisfies }$\lnot \lnot =$ 
\textit{id (see}$\left[ 13\right] $\textit{, page 270).}

\textit{Proof: } (i) $\Longleftrightarrow $ (ii) Suppose A $\subset $ X is
closed and p $\in $ A$^{C}$. \ Then p has a T$_{0}$- separation from each
point a $\in $ A. \ If X is pre-Hausdorff, then ($\forall $a $\in $ A)($%
\exists $ N$_{a}$, N$_{p_{a}}$ $\in $ $\tau $) such that a $\in $ N$_{a}$, p 
$\in $ N$_{p_{a}}$, and N$_{a}\cap $ N$_{p_{a}}$ = $\varnothing $. \ Then p $%
\in $ U = $\dbigcap\limits_{a\in A}$ N$_{p_{a}}$, A $\subset $ V = $%
\dbigcup\limits_{a\in A}$ N$_{a}$, and U $\cap $ V = $\varnothing $. \ Since
X is principal we have that U is open and, consequently, X is regular. \
Conversely, suppose that x, y $\in $ X have a T$_{0}$- separation by, say, U$%
_{x}\in $ $\tau $, where x $\in $ U$_{x}$ and y $\notin $ U$_{x}$. \ Then U$%
_{x}^{C}$ is closed so, if X is regular, there are disjoint open sets U and
V such that x $\in $ U and U$_{x}^{C}$ $\subset $ V. \ Thus, X is
pre-Hausdorff.

(ii) $\Longleftrightarrow $ (iii) See $\left[ 2\right] ,$ Theorem 2.9.

(iii) $\Longleftrightarrow $ (iv) In the topos of sheaves on X, the negation
operator $\lnot :\tau \rightarrow \tau $ is defined by $\lnot $ U = interior
(U$^{C}$) (see $\left[ 13\right] $, Chapter 2). \ Then $\lnot \lnot $ U =
interior($\overline{\text{U}}$), and so $\lnot \lnot $ U = U iff U =
interior($\overline{\text{U}});$ i.e., iff U is a regular open set (see $%
\left[ 4\right] $, page 92). \ It is easily shown that an open set is
regular iff it is clopen.

\bigskip

REMARKS 2.2.

\ \ (i) \ The proof of Theorem 2.1 shows that a regular space is
pre-Hausdorff even if it is not principal. \ Clearly the converse is false;
for if X is a Hausdorff space which is not T$_{3}$, then X is pre-Hausdorff
but not regular.

\ \ (ii) \ Also, a 0-dimensional space is pre-Hausdorff even if it is not
principal. \ However this is not true conversely; for the set of real
numbers \textbf{R }with the usual open interval topology is a (pre-)
Hausdorff space which is not 0-dimensional. \ In fact, dim(\textbf{R) }= 1
(see $\left[ 9\right] $, page 25, Example III).

\bigskip

In $\left[ 5\right] $ it is shown that a principal topological space (X, $%
\tau $) is regular if and only if the minimal basis for $\tau $ forms a
partition of X. \ Consequently, we have the immediate

\bigskip

COROLLARY 2.3. \textit{If X is a finite set, then the distinct pre-Hausdorff
topologies on X are in one-to-one correspondence with the distinct
partitions on X.}

\bigskip

In $\left[ 5\right] $ there is an algorithm using matrices, and a computer
program, to compute the number of regular (hence, pre-Hausdorff) topologies
on a finite set. \ Alternatively, several methods for counting the number of
partitions on a set with n-elements, the so-called "n-th Bell Number" B(n),
are well-known (see $\left[ 17\right] ,$, page 33). \ The 14th Bell number,
for instance, is B(14) = 190,899,322, which is accordingly the number of
distinct pre-Hausdorff topologies on a 14-element set as mentioned in the
introduction.

Finite (and other) pre-Hausdorff spaces can also be described using the
notion of a Borel field.

\bigskip

DEFINITION 2.4. \ A Borel field F (on a fixed set B) is a non-empty family
of subsets of B such that F is closed with respect to complements and
countable unions; i.e., F satisfies:

\ \ (i) if A $\in $ F, then A$^{C}$ $\in $ F

\ \ (ii) if $\left\{ \text{A}_{i}\right\} _{i=1}^{\infty }$ $\subset $ F,
then $\dbigcup\limits_{i=1}^{\infty }$ A$_{i}$ $\in $ F.

\bigskip

REMARKS 2.5.

\ \ (i) A Borel field is also known as a $\sigma $-algebra (see $\left[ 16%
\right] $, page 17, for instance).

\ \ (ii) If F is a Borel field on B, then clearly B $\in $ F and $%
\varnothing $ $\in $ F..

\ \ (iii) It follows immediately from DeMorgan's laws that a Borel field is
also closed with respect to countable intersections; i.e., if $\left\{ \text{%
A}_{i}\right\} _{i=1}^{\infty }$ $\subset $ F, then $\dbigcap\limits_{i=1}^{%
\infty }$ A$_{i}$ $\in $ F.

\ \ (iv) If F is a Borel field on B and F is countable, then (B, F) is a
topological space which has the following properties:

\ \ \ \ \ \ \ \ \ (1) arbitrary intersections of open sets are open, and

\ \ \ \ \ \ \ \ \ (2) every open set is also closed; i.e., every open set is
clopen (both open and closed).

\bigskip

COROLLARY 2.6. \ \textit{Suppose X is a finite set and }$\tau $ \textit{is a
family of subsets of X. \ }$\tau $ \textit{is a Borel field if and only if
(X, }$\tau $\textit{) is a pre-Hausdorff space.}

\textit{Proof: \ }Follows immediately from Remark 2.5 (iv) and Theorem 2.1.
\ Clearly this result is also true for any set X if $\tau $ is countable.

\bigskip

In $\left[ 20\right] $, Szekeres and Binet prove that the set of all Borel
fields on a finite set is in one-to-one correspondence with the number of
equivalence relations on that set. \ It is well known that the number of
equivalence relations on a finite set are in one-to-one correspondence with
the number of partitions on that set. \ Consequently, Corollary 2.6 is
equivalent to Corollary 2.3.

Of the 190,899,322 distinct pre-Hausdorff topologies on a set with 14
elements there are, of course, many which are homeomorphic. \ To
characterize homeomorphic pairs of finite pre-Hausdorff spaces, we look at
the basis consisting of the "minimal" open sets which, by the proof of
Corollary 2.3, forms a partition. \ For a finite pre-Hausdorff space X, we
shall denote this partition which generates X by P$_{\text{X}}$. \ The
following shows that for finite pre-Hausdorff spaces to be homeomorphic,
their generating partitions must "look" the same.

\bigskip

PROPOSITION 2.7. \ \textit{Finite pre-Hausdorff spaces (X, }$\tau $\textit{)
and (Y, }$\sigma $\textit{) are homeomorphic if and only if there exists a
bijective correspondence between P}$_{X}$ \textit{and P}$_{Y}$ \textit{that
preserves the cardinality of the corresponding blocks.}

\textit{Proof:} \ If X and Y are homeomorphic, then the condition on their
generating partitions follows immediately since a homeomorphism is a
bijective open mapping. \ Conversely, suppose there exists a bijective
correspondence between P$_{\text{X}}$ and P$_{\text{Y}}$ which preserves the
cardinality of the corresponding blocks, and that P$_{\text{X}}$ = $\left\{
B_{i}\right\} _{i=1}^{k}$ and P$_{\text{Y}}$ = $\left\{ C_{i}\right\}
_{i=1}^{k}$\textit{\ \ }are labeled so that B$_{i}$ and C$_{i}$ each have
the same cardinality for all $i$ = 1, 2, ..., $k$. \ Then, for each $i$, we
can choose a bijection f$_{i}$: B$_{i}$ $\rightarrow $ C$_{i}$. \ The
function f: X $\rightarrow $ Y defined by f(x) = f$_{i}$(x), \ for x $\in $ B%
$_{i}$, is clearly a homeomorphism.

\bigskip

It follows from Proposition 2.7 that the number of non-homeomorphic
pre-Hausdorff spaces on a set with n-elements is p(n) = the number of
partitions of n according to the following.

\bigskip

DEFINITION 2.8. \ A partition of a positive integer n is a finite
nonincreasing sequence of positive integers $\lambda _{1},$ $\lambda _{2},$
..., $\lambda _{r}$ such that $\dsum\limits_{i=1}^{r}\lambda _{i}$ = n.

\bigskip

The problem of computing p(n) in general is complex and has recieved much
attention from mathematicians, especially after the landmark paper by G. H.
Hardy and S. Ramanujan in 1918 ($\left[ 6\right] $). \ A comprehensive
summary of results can be found in $\left[ 1\right] $, where there is also a
table of values for p(n) up to p(100). \ Thus, the number of
non-homeomorphic pre-Hausdorff topologies on a set with 14 elements is p(14)
= 135 (see $\left[ 1\right] $, page 238).

\bigskip

Suppose X is a topological space, and B $\subset $ X. \ Recall that a point
b $\in $ X is called a generic point of B provided $\overline{\left\{ \text{b%
}\right\} }$ = B, and that X is called sober provided every closed
irreducible (i.e., cannot be decomposed into a union of two or more smaller
closed subsets) subset of X has a unique generic point. \ See $\left[ 10%
\right] ,$ page 230 for an interesting explanation of the term sober.

\bigskip

THEOREM 2.9. \ \textit{A topological space }(X, $\tau $) \textit{is
Hausdorff if and only if }X \textit{is both pre-Hausdorff and sober.}

\textit{Proof:}\ The implication to the right is immediate because Hausdorff
spaces are naturally pre-Hausdorff and, furthermore, they are always sober
(see $\left[ 13\right] ,$ page 475).

Now suppose that X is both pre-Hausdorff and sober, but not Hausdorff. \
Then $\exists $ x, y $\in $ X such that x $\neq $ y, and x and y have no T$%
_{2}$ separation in $\tau .$ \ Then x and y have no T$_{0}$ separation in $%
\tau $ either, which implies $\overline{\left\{ \text{x}\right\} }$ = $%
\overline{\left\{ \text{y}\right\} }.$ \ But then $\overline{\left\{ \text{x}%
\right\} }$ is a closed irreducible subset of X with more than one generic
point.

\bigskip

It is well-known that a topological space X is Hausdorff if and only if the
diagonal $\Delta _{\text{X}}$ (= $\left\{ \text{(x,x) : x}\in \text{ X}%
\right\} $) is closed in the product space X$^{2}$. \ Analogous results for
a pre-Hausdorff space (X,$\tau $) are given in terms of the following
relation R$_{0\text{ }}$on X.

\[
R_{0}\text{ = }\left\{ \text{(x, y) : x and y have no T}_{0}\text{
separation in }\tau \right\} \text{ }\subset \text{ X}^{2}
\]

Clearly R$_{0}$ is an equivalence relation.

\bigskip

THEOREM 2.10. \ \textit{The following are equivalent.}

\ \ \ (i) \textit{X is pre-Hausdorff.}

\ \ \ (ii) R$_{0}$ \textit{is closed in }$X^{2}$.

\ \ \ (iii) R$_{0}$ = $\overline{\Delta _{\text{X}}}$ \textit{.}

\ \ \ (iv) \ \textit{The quotient space }$\frac{\text{X}}{\text{R}_{0}}$ 
\textit{is Hausdorff}.

\textit{Proof:} \ We show that each of (ii), (iii), and (iv) is equivalent
to (i).

(ii) \ If R$_{0}$ is closed and x, y $\in $ \ X have a T$_{0}$ separation in 
$\tau $, then (x, y) $\notin $ R$_{0}$ = $\overline{\text{R}_{0}}$. \ So $%
\exists $ U, V $\in $ $\tau $ such that (x, y) $\in $ U$\times $V and (U$%
\times $V) $\cap $ R$_{0}$ = $\varnothing $. \ But this implies that U and V
are disjoint, for if p $\in $ U $\cap $ V, then (p, p) $\in $ (U$\times $V) $%
\cap $ R$_{0}$. \ Therefore x and y have a T$_{2}$ separation, and X is
pre-Hausdorff. \ Conversely, if X is pre-Hausdorff and (x, y) $\in $ R$%
_{0}^{C}$, then x and y have a T$_{2}$ separation in $\tau $ by , say, N$_{%
\text{x}}$ and N$_{\text{y}}$. \ Then (x, y) $\in $ N$_{\text{x}}\times $ N$%
_{\text{y}}$ $\subset $ R$_{0}^{C}$, so R$_{0}$ is closed.

(iii) \ That (iii) implies (i) follows immediately from (ii). \ For the
reverse implication, we have $\Delta _{\text{X}}$ $\subset $ \ R$_{0}$ $=$ $%
\overline{\text{R}_{0}}$\ (since X is pre-Hausdorff) so that $\overline{%
\Delta _{\text{X}}}$ $\subset $ R$_{0}$. \ For the reverse inclusion,
suppose a point (x, y) $\in $ R$_{0}$ = $\overline{\text{R}_{0}}$ has nbhd. U%
$\times $V in X$^{2}$. \ Then (x, x) $\in $ U$\times $V $\cap $ $\Delta _{%
\text{X}}$, so that (x, y) $\in $ $\overline{\Delta _{\text{X}}}$.

(iv) \ Suppose $\frac{\text{X}}{\text{R}_{0}}$ is Hausdorff, \ and that
distinct points x, y $\in $ X have a T$_{0}$ separation by a nbhd. of x. \
Then $\left[ \text{x}\right] $ (=$\left\{ \text{z: z R}_{0}\text{x}\right\} $%
) $\neq $ $\left[ \text{y}\right] $, so that $\left[ \text{x}\right] $ and $%
\left[ \text{y}\right] $ have a T$_{2}$ separation in $\frac{\text{X}}{\text{%
R}_{0}}$. \ Since the cannonical map q: X $\rightarrow $ $\frac{\text{X}}{%
\text{R}_{0}}$ is continuous, x and y have a T$_{2}$ separation in X. \
Conversely if X is pre-Hausdorff, then R$_{0}$ is closed by (ii). \
Furthermore, q: X $\rightarrow $ $\frac{\text{X}}{\text{R}_{0}}$ is easily
seen to be an open map. \ Consequently $\frac{\text{X}}{\text{R}_{0}}$ is
Hausdorff (see $\left[ 4\right] ,$ 1.6, page 140).

\bigskip

We now show that any space can be universally retracted onto a Hausdorff
space in the sense of adjunction as follows.

\bigskip

LEMMA 2.11. \ If (X, $\tau )$ is any topological space, (Y, $\sigma $) is a T%
$_{0}$ space, and f: X $\rightarrow $ Y is continuous, then f factors
uniquely through the quotient map q: X $\rightarrow $ $\frac{\text{X}}{\text{%
R}_{0}}$; i.e., $\exists !$ continuous $\overline{f}$: $\frac{\text{X}}{%
\text{R}_{0}}$ $\rightarrow $ Y such that f = $\overline{f}$ $\circ $ q.

\textit{Proof:} \ Define $\overline{f}$($\left[ \text{x}\right] $) = f(x). \
Then $\overline{f}$ is well-defined since Y is T$_{0}$, and $\overline{f}$
is continuous since $\frac{\text{X}}{\text{R}_{0}}$ is equipped with the
coinduced topology; i.e., the quotient topology in TOP.

\bigskip

THEOREM 2.12. \ \textit{The inclusion functor inc}$_{2,2}$\textit{: T}$_{2}$%
\textit{-Top }$\hookrightarrow $ \textit{T}$_{0,2}$\textit{-Top has a left
adjoint L}$_{2,2}$ \textit{which is a retract.}

\textit{Proof: \ }Define L$_{2,2}$(X) = $\frac{\text{X}}{\text{R}_{0}}$. \
By Lemma 2.11, the quotient map q: X $\rightarrow $ $\frac{\text{X}}{\text{R}%
_{0}}$ provides a universal arrow from any pre-Hausdorff space X to the
Hausdorff space $\frac{\text{X}}{\text{R}_{0}}$. \ The object $\frac{\text{X}%
}{\text{R}_{0}}$ and the universal arrow q completely determine the left
adjoint to $inc_{2,2}$ (see $\left[ 12\right] $, Theorem 2 (ii), page 81). \
L$_{2,2}$ is a retract by Theorem 2.10 (iii).

\bigskip

COROLLARY 2.13. \ \textit{The inclusion functor inc}$_{2}$\textit{: T}$_{2}$%
\textit{\ - Top }$\hookrightarrow $ \textit{Top has a left adjoint L}$_{2}$ 
\textit{which is a retract.}

\textit{Proof: \ }Combining Corollary 1.7 with Theorem 2.12, we\textit{\ }%
define L$_{2}$ = L$_{2,2}\circ $ L$_{0,2}$.

\bigskip

The functor L$_{2}$ can be described without the use of L$_{0,2}$ and L$%
_{2,2}$. \ To this end, we now construct L$_{2}$ directly by way of forming
quotients by an equivalence relation. \ We take a general approach which
also shows T$_{0}$ -TOP and T$_{1}$ - TOP to be reflective, and gives an
explicit description of the left adjoints to their inclusions into TOP.

\bigskip

DEFINITION 2.14. \ Let (X, $\tau $) be a topological space. \ For each $i$ =
0, 1, 2, define a relation R$_{i}$ on X by:

\[
(\text{x, y) }\in \text{ R}_{i}\text{ iff }\forall \text{Y }\in \text{ T}_{i}%
\text{ - TOP, }\forall \text{ continuous f : X }\rightarrow \text{ Y, f(x) =
f(y).}
\]

\bigskip

REMARK 2.15. \ R$_{0}$ as defined in 2.14 equals R$_{0}$ as defined above.

\bigskip

LEMMA 2.16. \ \textit{For each i = 0, 1, 2 we have the following:}

\ \ \ \ (i) \ \textit{R}$_{i}$\textit{\ is an equivalence relation.}

\ \ \ \ (ii) \ \textit{If Y }$\in $ \textit{T}$_{i}$ - \textit{TOP and f : X 
}$\rightarrow $ \textit{Y is continuous, then f factors through the quotient
map q : X }$\rightarrow $ $\frac{X}{R_{i}}$ .

\ \ \ \ (iii) \ $\frac{X}{R_{i}}$ $\in $ \textit{T}$_{i}$ - \textit{TOP.}

\textit{Proof: \ }(i) Straightforward.

\ \ \ (ii) \ Given a continuous function f : X $\rightarrow $ Y with Y $\in $
T$_{i}$ - TOP, define $\overline{\text{f}}$ : $\frac{\text{X}}{\text{R}_{i}}$
$\rightarrow $ Y by $\overline{\text{f}}$($\left[ \text{x}\right] $) = f(x).
\ Then $\overline{\text{f}}$ is well-defined by definition of R$_{i}$, and $%
\overline{\text{f}}$ is continuous since $\frac{\text{X}}{\text{R}_{i}}$ has
the quotient topology in TOP.

\ \ \ (iii) \ Suppose that $\left[ \text{x}\right] $ $\neq $ $\left[ \text{y}%
\right] $ in $\frac{\text{X}}{\text{R}_{i}}$. \ Then $\exists $ Y $\in $ T$%
_{i}$ - TOP and $\exists $ continuous f : X $\rightarrow $ Y with f(x) $\neq 
$ f(y), which implies that f(x) and f(y) have a T$_{i}$ - separation in Y. \
Then $\left[ \text{x}\right] $ and $\left[ \text{y}\right] $ have a T$_{i}$
- separation in $\frac{\text{X}}{\text{R}_{i}}$ via inverse image of $%
\overline{\text{f}}$ : $\frac{\text{X}}{\text{R}_{i}}$ $\rightarrow $ Y.

\bigskip

THEOREM 2.17. \ \textit{For each i = 0,1,2, the inclusion functor inc}$_{i}$ 
\textit{: T}$_{i}$ \textit{- TOP }$\hookrightarrow $ \textit{TOP has a left
adjoint L}$_{i}$ \textit{: TOP }$\rightarrow $ \textit{T}$_{i}$ \textit{-
TOP. \ Moreover, each L}$_{i}$ \textit{is a retract.}

\textit{Proof: \ }Define L$_{i}$((X, $\tau $)) = $\frac{\text{X}}{\text{R}%
_{i}}$\bigskip . \ Then, by Lemma 2.16, L$_{i}($(X, $\tau $)) $\in $ T$_{i}$
-TOP, and the quotient map q : X $\rightarrow $ $\frac{\text{X}}{\text{R}_{i}%
}$ is universal among all arrows from X into a T$_{i}$ - space. \ If X $\in $
T$_{i}$ - TOP then, clearly, L$_{i}$(X) = X.

\bigskip

The functor L$_{0,2}$ of Corollary 1.7 can also be explicitly described
using the equivalence relation R$_{2}$. \ Indeed, if (X, $\tau $) is a
topological space, then $\frac{\text{X}}{\text{R}_{2}}$ is Hausdorff by 2.16
(iii). \ So (X, $\tau _{2}$ ), the topological space induced from $\frac{%
\text{X}}{\text{R}_{2}}$ via q : X $\rightarrow $ $\frac{\text{X}}{\text{R}%
_{2}}$ will be pre-Hausdorff. \ It is readily shown that the assignment (X, $%
\tau $) $\longmapsto $ (X, $\tau _{2}$) is left adjoint to the inclusion T$%
_{0,2}$ - TOP $\hookrightarrow $ TOP.

\bigskip

\textbf{3. Replacing Hausdorff With Pre-Hausdorff}

\bigskip

There are many known results in topology which concern Hausdorff spaces. \
Given such a result, a natural question is whether or not the result remains
true when Hausdorff is replaced with pre-Hausdorff. \ In this section we
point out some standard Theorems which can be generalized to the
pre-Hausdorff setting, and some which cannnot.

Pre-Hausdorff topologies share some invariance properties with Hausdorff
topologies:

\bigskip

PROPOSITION 3.1. (i) \textit{Each subspace of a pre-Hausdorff space is also
pre-Hausdorff.}

(ii) \textit{The Cartesian product of pre-Hausdorff spaces is also
pre-Hausdorff.}

\textit{Proof: }The proof of (i) is straightforward, (ii) follows
immediately from Theorem 1.6.

\bigskip

In the following result, as in $\left[ 11\right] $ page 112, we make a
distinction between normal spaces and T$_{4}$ spaces; namely, normal spaces
need not have closed points. \ Recall that every compact Hausdorff space is
normal (see $\left[ 11\right] $, p141).

\bigskip

THEOREM 3.3. \ Every compact pre-Hausdorff space is normal.

\textit{Proof: \ }Suppose\textit{\ }(X\textit{, }$\tau $) is a compact
pre-Hausdorff space. \ Since the Theorem is trivially true when $\tau $ is
the indiscrete topology, we assume that it is not and choose a closed set A $%
\subset $ X and a point x $\notin $ A. \ For each y $\in $ A, A$^{C}$
provides a T$_{0}$ separation of x and y. \ Since X is pre-Hausdorff, x and
y have a T$_{2}$ separation by , say N$_{\text{y}}$ $\ni $ y and N$_{\text{%
x,y}}$ $\ni $ x. \ Then $\left\{ \text{N}_{\text{y}}\right\} _{\text{y}\in 
\text{A}}$ is an open cover of A. \ Since A is compact, $\exists $ y$_{1}$,
... , y$_{\text{n}}$ $\in $ A such that U = $\dbigcap\limits_{i=1}^{n}$ N$_{%
\text{x,y}_{i}}$ and V = $\dbigcup\limits_{i=1}^{n}$ are open, and they
provide a disjoint separation of x and A.

Now suppose that A and B are disjoint, closed sets in X. \ By the above we
have, $\forall $ a $\in $ A, $\exists $ U$_{\text{a}}$, U$_{\text{a,B}}$ $%
\in $ $\tau $ such that a $\in $ U$_{\text{a}}$, B $\subset $ U$_{\text{a,B}%
} $, and U$_{\text{a}}$ $\cap $ U$_{\text{a,B}}$ = $\varnothing $. \ Since A
is compact, $\exists $ a$_{1}$, ... , a$_{n}$ such that U$_{\text{A}}$ = $%
\dbigcup\limits_{i=1}^{n}$ U$_{\text{a}_{i}}$ and V$_{\text{B}}$ = $%
\dbigcap\limits_{i=1}^{n}$ U$_{\text{a}_{i}\text{,B}}$ are both open, and
they provide a disjoint separation of A and B.

\bigskip

A special property of Hausdorff topologies is the following:

\[
\text{Each finite subset of a Hausdorff space is closed.}
\]

This is clearly not shared with pre-Hausdorff topologies; for if X is a
finite indiscrete space with more than one element, for instance, and A is
any non-void proper subset of X, then A is not closed.

\bigskip

Many important results involve mating compactness with the Hausdorff
property. \ An intriguing feature of compact Hausdorff spaces is that they
are essentially algebraic. \ Indeed, a well-known result is that the
category of compact Hausdorff spaces is algebraic (or "monadic") over the
category of sets (see $\left[ 12\right] $, chapter 6). \ A fact which is
crucial in proving this is the following:

\ \ \ \ If X is compact and Y is Hausdorff, then any continuous function f :
X $\rightarrow $ Y is a closed map.

This result is clearly false for pre-Hausdorff spaces; for example, if we
map a compact space (X, $\tau $) which is not indiscrete into X with the
indiscrete topology by the identity function, then we have a continuous
bijection which is not a closed map. \ Consequently, this identity map is
not a homeomorphism.

\ It is easily shown that if a category A is algebraic over a category B
(i.e., A is isomorphic to a category of T-algebras, where T is a monad in X
determined by an adjunction), and L $\dashv $ R : A $\rightarrow $ B is the
adjoint pair of functors which determines the isomorphism, then A satisfies:
if f : a $\rightarrow $ b is any morphism in A, and R(f) : R(a) $\rightarrow 
$ R(b) is an isomorphism in B, then f is an isomorphism in A. \ In the case
of compact Hausdorff spaces over the category of sets, this reveals the
well-known fact that a continuous bijection of compact Hausdorff spaces is a
homeomorphism. \ Since the example in the preceeding paragraph shows a
continuous bijection from a compact space to a pre-Hausdorff space which is
not a homeomorphism, we conclude that the category of compact pre-Hausdorff
spaces is not algebraic over the category of sets.

\bigskip

\textbf{References}

\bigskip

1. \ Andrews, G. E.: \textit{The Theory of Partitions, }Cambridge University
Press, 1998.

2. \ Arenas, F. G.: Alexandroff Spaces, \textit{Acta Math. Univ. Comenianae} 
\textbf{68} (1999), 17-25.

3. \ Das, P. and Banerjee, A. K.: Pairwise Borel and Baire Measures in
Bispaces, \textit{Arc. Math. (Brno)} \textbf{41} (2005), 5-15.

4. \ Dugundji, J.:\textit{\ Topology,} Allyn and Bacon, 1966.

5. \ Farrag, A. S. and Abbas, S. E.: Computer Programming for Construction
and Enumeration of all Regular Topologies and Equivalence Relations on
Finite Sets, \textit{Appl. Math. Comp. }\textbf{165 }(2005), 177-184.

6. \ Hardy, G. H. and Ramanujan, S.: Asymptotic formulae in Combinatory
Analysis, \textit{Proc. London Math. Soc.} (2) \textbf{17}, 75-115.

7. \ Herrlich, H.: Compact T$_{0}$ Spaces, \textit{Applied Categorical
Structures }\textbf{1} (1993), 111-132.

8. \ Herrlich, H.: A Concept of Nearness, \textit{Gen. Top. \& Appl.}\textbf{%
\ 5} (1974), 191-212.

9. \ Hurewicz, W. and Wallman, H.: \textit{Dimension Theory, }Princeton U.
Press, 1974.

\noindent\ \ \ \ 10. \ Johnstone, P.T.: \ \textit{Topos Theory,} L.M.S.
Math. Monograph No. 10, Academic Press, 1977.

\noindent\ \ \ \ 11. \ Kelley, J. L.: \textit{General Topology},
Springer-Verlag, 1985.

\noindent\ \ \ \ 12. \ MacLane,S.: \textit{Categories for the Working
Mathematician,} Springer-Verlag, 1971.

\noindent\ \ \ \ 13. \ MacLane, S. and Moerdijk, I.: \textit{Sheaves in
Geometry and Logic, }Springer-Verlag, 1992.

\noindent\ \ \ \ 14. \ Mielke, M. V.: Separation Axioms and Geometric
Realizations, \textit{Indian J. Pure Appl. Math. }\textbf{25} (1994),
711-722.

\noindent\ \ \ \ 15. \ Preuss, G.: \textit{Theory of Topological Structures,}
D. Reidel Publishing Co., 1987.

\noindent\ \ \ \ 16. \ Royden, H. L.: \textit{Real Analysis,} Macmillan
Publishing Co., Inc., 1968.

\noindent\ \ \ \ 17. \ Stanley, R. P.: \ \textit{Enumerative Combinatorics,}
Vol. 1, Wadsworth \& Brooks, 1986.

\noindent\ \ \ \ 18. \ Steiner, A. K.: The Lattice of Topologies: Structure
and Complementation, \textit{Trans. Am. Math. Soc. }\textbf{122} (1966),
379-398.

\noindent\ \ \ \ 19. \ Stine, J.: \ \textit{Pre-Hausdorff Objects in
Topological Categories,} Ph. D. Dissertation, University of Miami, 1997.

\noindent\ \ \ \ 20. \ Szekeres, G. and Binet, F.E.: Notes On Borel Fields
Over Finite Sets, \textit{Ann. Math. Stat. }\textbf{29 }(1957), 494-498.

\end{document}